\input amstex
\documentstyle{amsppt}

\define\RR{\Bbb R}
\define\CC{\Bbb C}
\define\SS{\Bbb S}
\define\PP{\Bbb P}

\topmatter
\title Real Polynomial Rings and Domain Invariance
\endtitle
\rightheadtext {}
\author Jon A. Sjogren
\endauthor
\affil Inphormax LLC
\endaffil
\address Sullivan Station, Somerville MA
\endaddress
\date 22 December 2014
\enddate
\endtopmatter
\document

\abovedisplayskip=13pt
\belowdisplayskip=13pt

\head Introduction \endhead

Among the different types of proof for the ``Fundamental Theorem of Algebra'' (FTA or Main Theorem on Real/Complex Polynomials), some controversy has arisen as to the applicability of Brouwer's celebrated Fixed-Point Theorem (FPT). On the one hand D. Reem in \cite{7} points out that efforts of the past (see \cite{5}) have been revealed as inadequate, and it has been argued in \cite{6} that the endeavor is actually futile.

On the other hand, it is well-understood that a purely algebraic proof of FTA should not exist. One must at least use an analytic or topological step in the demonstration; one way is to demonstrate that an odd-degree real polynomial possesses a real root. Actually, Brouwer's FPT would seem to fulfill admirably this need for an ``analytical ingredient'' in the proof. In fact, both  \cite{4} and \cite{7} assert that such a proof might be carried out, despite the claims of ``impossibility''.

We maintain that one can indeed prove FTA using Brouwer's FPT. The road is a bit indirect, meandering through Brouwer's other pioneering result on Euclidean topology, Invariance of Domain.
You have proved FTA when you prove that any irreducible {\it real} polynomial has degree one or two.
The existence of an {\it irreducible} polynomial of higher degree $n$ leads to a real vector space of that same dimension on which the polynomial product induces a ring structure which in fact yields an integral domain.

Consider the squaring map on this real algebra, and scale it to map the entire vector space to a Cartesian sphere of dimension one less.  Scaling the domain vector variable also, we get a squaring mapping from unit sphere to unit sphere, in fact from the real projective space $\PP^{n-1}$ to the sphere 
$\SS^{n-1}$. The mapping is topologically continuous (on Hausdorff spaces) due to the bilinearity of the original product operation.

The integrity condition (entire domain) shows that this mapping is injective. All spaces in sight are compact Hausdorff, so such a 1-1 mapping induces a homeomorphism onto the image.
If one throws in the ``connectedness'' of the spaces involved, then the (Brouwer) Invariance of Domain Theorem implies that in fact the image of $\PP^{n-1}$ must be the whole sphere.  The inference that these two spaces are homeomorphic leads to the conclusion that $n \leq 2$.

To show that $\PP^{n-1}$ is not topologically the same as $\SS^{n-1}$, construct a non-trivial loop in $\PP^{n-1}$
by projecting a longitude running South pole to North pole, from its universal covering sphere.
This loop could not be homotopic to the constant loop (at the chosen base point) since the homotopy would lift to the Universal Cover keeping end-points invariant.

What does this have to do with Brouwer's FPT?  As Prof. Terence Tao points out on his discussion page \cite{10}, the latter theorem is necessary in some form to prove Invariance of Domain.

\head  Self-Mapping on the Space of Real Polynomials\endhead

The last decade or so has seen remarkable simplifications in the proofs of classic theorems of algebra. Not only are the new proofs easier to understand, they lend themselves better to symbolic computation and to applications. For further insight, one may consult \cite{11} and \cite{12} among other sources. Consider for instance that pillar of algebraic geometry, Hilbert's {\it Nullstellensatz}, or Theorem of Zero Placement. Since the basic result applies to equations defined over an algebraically closed field, it is at least aesthetic to consider proofs of FTA (algebraic closure of the complex numbers $\CC$) that are concordant with some of the new proofs of the {\it Nullstellensatz}, such as that of Almira in \cite{1}. The present note also owes to observations by J. Vonk in \cite{2} and E. Arrondo in \cite{12}.

Considering the variety of approaches to FTA, it is of interest to provide a proof of FTA that employs  unexpected but not especially advanced tools of topology and functional analysis.
We formulate the problem as, given $p(x) \in \RR[x]$ of degree $n > 2$, to show that $p$ is reducible (equals the product of factors each of degree $1 \leq      k \leq n-1$). This factorizability is considered obvious when $n$ is an odd integer, so we may assume $n$ even. The structure of the argument is the same either way.

We are actually finished with the proof as soon as we show that an {\it irreducible} $p(x) \in \RR[x]$ must have degree $n \leq 2$. Starting with such a polynomial (but non-constant) we may work with a quotient ring
$$\Cal O:= \RR[x]/p{(x)}$$,
which under ``polynomial addition, scalar and polynomial multiplication'' is an integral domain (entire commutative ring), in fact an $\RR$-algebra. Regarding the terms $X^k,\,\, k = 1, \dotsc, n-1$ as a linear basis, the structure of $\Cal O$ leads to a product
$$\circ : \RR^n \times \RR^n \rightarrow \RR^n$$
with canonical unit element and no zero divisors. This (polynomial) product is linear in the first and second vector variable, so bi-linear, and with the metric topology inherited from $\RR$, it is continuous in both $n$-vector variables.

For $a(x) \in \Cal O$ the mapping $f(a) = a \circ a$ is well-defined and continuous on $\RR^n$. Restricting $f$ to vectors of norm $= 1$, we arrive at a mapping $\psi = \SS^{n-1} \rightarrow \SS^{n-1}$ defined by
$$\psi(a) = \frac{a \circ a}{|a \circ a|}.$$
Here $|a \circ a|$ really means
$$\left\{(a \circ a) \cdot (a \circ a)\right\}^{1/2},$$
using the Cartesian inner product on $\RR^n$.

We summarize the algebra which has been used up to this point. An irreducible $p(x)$ generates a principal, prime ideal $(p)$. The coordinate (quotient) ring $\Cal O = \RR[x]/\left(p(x)\right)$ is entire, since any zero divisor would be a factor of $p(x)$ that is neither a scalar nor a scalar multiple of $p(x)$. The ring $\Cal O$ is also an $\RR$-vector space with basis the cosets
$$e_0 = \left\{1 + (p(x))\right\}, \, e_1 = \left\{x + (p(x))\right\}, \dotsc , e_{n-1} = \left\{x^{n-1} + (p(x))\right\}.$$
Consequently, $\Cal O$ is an $\RR$-algebra, and in fact, a {\it division algebra} (where any non-zero element $b$ has an multiplicative inverse $c = b^{-1}$ where $b \circ c = c \circ b = \bold 1$, the coset containing the scalar $1 \in \RR$). However, it is worth remarking that our results follow using only the weaker condition of ``no zero divisors'' in place of ``all non-zero elements are invertible''.

Next we show through topology how such an algebraic structure imposes a restriction on the irreducible $p(x)$, namely the size of its degree. Starting with the continuous mapping $\psi : \SS^{n-1} \rightarrow \SS^{n-1}$ we note that this induces $\hat{\psi} = \PP^{n-1} \rightarrow \SS^{n-1}$ on {\it real projective} $(n-1)$-space, in view of the observed fact $\psi(y) = \psi (-y)$ for all $y \in \SS^{n-1}$.

The resulting mapping $\hat{\psi}$ turns out to be {\it injective}, since if $\hat{\psi}(u) = \hat{\psi}(v)$, using ``liftings'' $u, v \in \SS^{n-1}$, we infer
$$|v^2| u^2 = |u^2| v^2$$,
hence  $u^2 - \alpha^2 v^2 = 0$ where
$$\alpha^2 = |u \circ u|/|v \circ v|.$$
Since $\Cal O$ is an integral domain, $0 = (u - \alpha v) \circ (u + \alpha v)$ implies $u = \alpha v$ or $u = -\alpha v$.
But $|u| = |v| = 1$ so $\alpha = \pm 1$ and $u = v$ as elements of $\PP^{n-1}$.

Therefore
$$\hat{\psi} : \PP^{n-1} \rightarrow \SS^{n-1}$$
is a continuous, one-to-one mapping of compact Hausdorff spaces. Under the conditions $f: Y \to Z$, where $f$ is a (continuous) mapping, $Y$ is compact (Hausdorff) and $Z$ also has the  Hausdorff property, we immediately conclude that $f$ is a {\it closed} map. But by definition, a continuous and closed {\it bijective} function is a homeomorphism. Hence $\hat{\psi}$ is a homeomorphism onto its image $\psi(\PP^{n-1}) \subset \SS^{n-1}$.
For derivations of basic topological facts one may consult e.g. \cite{9}, p. 226.

We revert to $f: Y \to Z$, where $f$ is $1-1$ (an embedding) but in this case both $Y$ and $Z$ are compact $m$-manifolds, and $Z$ is a connected space.
The conclusion (see Hatcher, \cite{3}, p. 172) is that $f$ is also surjective, hence a {\it homeomorphism}.

The proof of this general remark rests on the Theorem of Invariance of Domain, due to Brouwer. The image $f(Y)$ is closed in $Z$ since $f(Y)$ is compact and $Z$ is Hausdorff. Now we establish that also $f(Y)$ is {\it open} in $Z$. Given $y \in Y$ choose a neighborhood $V \subset Z$ of $f(y)$ that is homeomorphic to $\RR^m$, together with a neighborhood $U \subset Y$ of $y$ with $U \simeq \RR^m$ and small enough so that $f(U) \subset V$. Brouwer's theorem states that $f(U)$ is open in $V$, hence in $Z$, so an arbitrary       $f(y)$ has an open neighborhood of $Z$, entirely contained in $f(Y)$.
Hence $f(Y)$ is open and must equal its connected component in $Z$.

Applying this general proposition to the case of the one-to-one mapping $$\hat{\psi}: \PP^{n-1} \to \SS^{n-1},$$ taking $m = n-1$, we conclude since the spaces are compact (Hausdorff) and $\SS^{n-1}$ is connected (for $n >1$), that $\hat{\psi}$ is surjective and a homeomorphism.
But for $n >2$, the two spaces possess unlike topological invariants. Therefore an irreducible polynomial $p(x)\in \RR[x]$ of degree $n \geq 3$ is impossible which was to be proved.

We conclude by showing that $\PP^{n-1} \simeq \SS^{n-1}$ implies $n \leq 2$ without employing higher topological machinery.
We can now form the composition of maps
$$\SS^{n-1}\overset\pi\to\longrightarrow \PP^{n-1} \overset\hat{\psi}\to\longrightarrow \SS^{n-1}$$
where $\pi$ is the {\it projection} (identification map $\{y, -y\}$ in the Universal Cover). In $\SS^{n-1}$ we have South and North poles, so let $\beta$ be the Prime Meridian arc connecting $\text{South}$ to $\text{North}$ in $\SS^{n-1}$ and $\gamma$ be the loop $\pi (\beta)$ in $\PP^{n-1}$. In other words $\gamma$ is a closed arc beginning and ending at $\{\text{South}, \text{North}\}$. Since whenever $n > 2$ all loops in $\SS^{n-1}$ are contractible within $\SS^{n-1}$ to a constant arc, we obtain by the use of homeomorphism $\hat{\psi}$, a homotopy $H_t$ on $\PP^{n-1}$ where $H_0$ is $\gamma$ and $H_1$ is the constant arc $H_1: [0,1] \to \{\pi(\text{South})\}$.

By the unique homotopy covering theorem with unique path lifting (see \cite{3}, p. 60), the homotopy $H: I^2 \to \PP^{n-1}$ has a lifting $\hat{H}$ to $\SS^{n-1}$ where the endpoints of $\hat{H}_t$ for varying $t$ remain fixed. But the endpoints of $\hat{H}_1$ are both $\{\text{South}\}$ and the endpoints of $\hat{H}_0$ are 
$\{\text{South}\}$ and $\{\text{North}\}$, from which contradiction we conclude that $n \leq 2$.

\bigskip
{\bf Remark.}  The discussion in \cite{10} shows how Brouwer's FPT seems to be pivotal in proving Brouwer's Theorem on the Invariance of Domain. Historically this represents a juncture in the above proof of the ``Fundamental Theorem of Algebra'' that seemed obvious, but in fact requires careful verification.

\enddocument